\documentclass[11pt,a4paper]{article}
\usepackage{amsmath}
\usepackage{graphicx}
\usepackage{amssymb}
\input amssymb.sty
\renewcommand{\le}{\leqslant}
\renewcommand{\ge}{\geqslant}

\def \N{{\mathbb  N}}
\def\C{{\mathbb C}}

\def \R{{\mathbb R}}

\def\E{\mathbb{E}}
\def\I{\mathcal{I}}
\def\J{\mathcal{J}}
\def\D{\mathcal{D}}

\def\asc{\mathrm{asc}}
\def\desc{\mathrm{desc}}

\def\crr{\mathrm{cr}}
\def\trr{\mathrm{tr}}

\def\oml{\omega_{\cal L}}
\def\dz{\frac{\rm d}{{\rm d}z}}
\def\ds{\frac{\rm d}{{\rm d}s}}

\def\A{{\cal A}}
\def\B{{\cal B}}
\def\bu{\bullet_\kappa}
\newtheorem{theorem}{Theorem}
\newtheorem{lemma}{Lemma}
\newtheorem{corollary}[lemma]{Corollary}
\newtheorem{definition}[lemma]{Definition}
\newtheorem{proposition}[lemma]{Proposition}

\def\endproof{\hfill $\square\qquad$\endtrivlist} 

\newtheorem{remark}[lemma]{Remark}
\newtheorem{example}[lemma]{Example}

\begin{document}
\title{Simply Generated Trees, B-Series\\
and Wigner Processes}

\author{ Christian Mazza\\
          Section de Math\'ematiques\\
	  2-4 Rue du Li\`evre, CP 240\\
	  CH-1211 Gen\`eve 24 (Suisse)\\
	   christian.mazza@math.unige.ch}
\date{}
\maketitle
\renewcommand{\thefootnote}{}

\begin{abstract}
We consider simply generated trees,
like rooted plane trees, and consider the problem
of computing generating functions of so-called bare
functionals, like the
tree factorial, using B-series from Butcher's theory.
 We exhibit
a special class of functionals from  probability theory: the associated
generating functions can be seen as limiting
traces of product of semi-circular elements.
\end{abstract}


\medskip\noindent
2000 AMS Subject Classification.
Primary 05a15 Secondary 05c05, 30b10, 46l54, 81s25
\section{Introduction\label{s.intro}}
Let ${\cal F}_n$ denote the set of rooted plane trees 
of size $n$. Simply generated trees are families of trees obtained by
assigning weights $\omega(t)$ to the elements $t\in {\cal F}=\cup_n {\cal F}_n$
using a degree function $\psi(z)=1+\sum_{k\ge 1}\psi_k z^k$
(see \cite{MeirB}). Basically,
the weight $\omega(t)$ of some $t\in {\cal F}$ is obtained by multiplying
the factors $\psi_{d(v)}$ over the nodes $v$ of $t$,
where $d(v)$ denotes the outdegree of  $v$.
Our main topic is the study of generating functions
$$Y(z)=\sum_{t\in {\cal F}}\omega(t)B(t)z^{\vert t\vert},$$
associated with multiplicative functions $B:{\cal F}\longrightarrow \R$
defined recursively by using a sequence of real numbers
$\{B_k\}_{k\in \N^+}$. We call such multiplicative functions
{\it bare Green functions}: $\sum_{t\in {\cal F}_n}B(t)\omega(t)$
represents the sum of the Feynman amptitudes associated to the relevant diagrams of size $n$
in some field theory, and the generating function is then
a part of the perturbative expansion of the solution of some
equation describing the system
(see \cite{Broadhurst,Brouder,KreimerB,Kreimer}
).

 In Section \ref{Runge}, we
 give   an equation satisfied by $Y$ when 
the weights $B_k$ come from some master function
$L(z)=\sum_{m\ge 0}L_m z^m$, with
$B_k\equiv L(k)/k$, $\forall k\in \N^+$. We use series indexed by trees, the so-called
B-series, as defined in \cite{Wanner,Hairer},
to show in Theorem \ref{DifferentialSys} that $Y$ solves
$$Y'=L(1+\theta)\Psi(Y),$$
 where
$\theta$ is the differential operator $\theta=z {\rm d}/{\rm d}z$.
 \cite{berg} considers a similar problem
for {\it additive tree functionals}
$s(t)$ defined on varieties of increasing trees, like
 $s(t)=\ln(B(t))$. Assuming  some constraints on
 the degree function $\Psi(z)$, it is proven that the exponential
 generating function 
 $$S(z)=\sum_{t\in {\cal F}}\omega(t)s(t)z^{\vert t\vert}/\vert t\vert !,$$
 is given
 by the formula 
$$S(z)=W'(z)\int_0^z (F'(u)/W'(u)){\rm d}u,$$
 where $F(u)=\sum_{m\ge 0}\ln(B_m)W_m u^m/m!$ and
 $W(z)=\sum_{m\ge 0}W_m z^m/m!$ solves $W'=\Psi(W)$. 
 We also consider a central functional called
 the {\em tree factorial}, denoted by $t!$  in
 the sequel, 
which is relevant in various fields,
like algorithmics \cite{Dob,Mah},
stochastics \cite{Fill,MeirA},
numerical analysis (see for example \cite{ButcherA,Hairer}),
and physics \cite{Brouder,Kreimer}. 
We focus on its negative powers $1/(t!)^{l+1}$, $l\in\N$, which
 do not admit a master function when $l\ge 1$. 
\cite{Brouder}
 solved the case $l=1$ by using the so-called Butcher's group
 (see for example \cite{Wanner,Hairer}). We  provide in Theorem
\ref{Inversion}
 a differential equation for the associated generating function,
 $\forall l\in\N$.

 In Section \ref{WignerProcesses},
 we define special multiplicative
 functionals for which the weights $B_k$ are related to the
 covariance function $r$ of some gaussian process, as
 $B_k=\beta^2 r(2k-1)$, for some positive constant $\beta>0$. We show that the generating
 function $Y$ is related to the mean normalized trace
 of products of large symmetric random matrices having
 independent and indentically distributed versions of the  process as entries.
 Theorem \ref{t.wig} gives then  a  differential equation for
 the evolution of the trace of a stationary
 Wigner processes.  
 It follows that most of the examples given in \cite{Broadhurst,Kreimer}
can be expressed in terms of traces of large
 random matrices.
 In Section \ref{Triangular},
 we show how B-series can be useful for studying 
 traces of triangular operators appearing in free probability.

\section{Basic notions\label{s.def}}
A rooted tree $t\in {\cal R}$ is a triple $t=(r,V,E)$ such that i) $(V,E)$ is a non-empty directed tree with node set
$V$ and edge set $E$, ii) all edges are directed away from the root $r\in V$.
The set of rooted trees of order $n$ is denoted by ${\cal R}_n$, and the
set of rooted trees is ${\cal R}=\cup_n {\cal R}_n$.
A rooted plane tree $t\in {\cal F}$ is a quadruple $t=(r,V,E,L)$ satisfying i) and ii) 
and iii) 
$L:=\{(\{w:\ vw\in E\},\ L_v):\ v\in V\}$ is a collection of $\vert V\vert$ linear orders. Given $v\in V$, let
${\rm ch}(v):=\{w:\ vw\in E\}$ be the set of children of $v$. $d(v):=\vert{\rm ch}(v)\vert$ is the  outdegree of $v$. A rooted planar tree
can be seen in the plane with the root in the lowest position, such that the orders $L_v$ coincide with
the left-right order. Next consider the partial ordering $(V,\le)$ defined by $u\le v$ if and only if
$u$ lies on the path linking $r$ and $v$. Given $v\in V$ and $t\in {\cal R }   $
let $t_v$ be the subtree
of $t$ rooted at $v$ spanned by the subset $\{w;\ v\le w\}$. 
A rooted labelled tree is a quadruple $t=(r,V,E,l)$ satisfying i) and ii), with a labelling 
$l:\ V\setminus \{r\}\longrightarrow [\vert V\vert ]:=\{1,\cdots,\vert V\vert\}$ such that
$l(u)<l(v)$ when $u<v$.  The set of rooted labelled trees of order $n$ is denoted by
${\cal L}_n$. Let ${\cal L}=\cup_n {\cal L}_n$. This family is a special variety of increasing trees, as
defined in \cite{berg,Fla}.

We next assign weights to the elements of ${\cal F}_n$, the set of rooted planar trees of order $n$:
the resulting family of trees is said to be simply generated
(see \cite{MeirA}).
Given a sequence $\psi=\{\psi_k\}_{k\in \N}$ of real numbers with
$\psi_0=1$, define recursively the weight $\omega(t)$ of $t\in{\cal F}$ as
$$\omega(t)=\psi_k \prod_{i=1}^k \omega(t_i),\ k=d(r),\ \omega(t)=\prod_{v\in V}\psi_{d(v)}.$$
where $t_1,\cdots,t_k$ are the $d(r)$ subtrees of $t$ rooted
at ${\rm ch}(r)$.
Let $\psi(z):=1+\sum_{k=1}^\infty \psi_k z^k$ be the generating function
of the weight sequence $\psi$. 
Our favourite example is
$\psi(z)=1/(1-z)$, with $\omega(t)=1$, $\forall t$
(see \cite{berg,MeirB} for various interesting choices).

We will be concerned with
 functionals $B:{\cal F}\longrightarrow \R$, where ${\cal F}=\cup_n {\cal F}_n$,
  called {\it bare Green functions}. This terminology is taken from
quantum field theory where bare Green functions occur during the action
of the renormalization group (see for example
\cite{Broadhurst}, \S~4.2 or
 \cite{Brouder} , \S~6.1).
Let ${\bf B}$ denote the set of bare Green functions.
Any element $B\in {\bf B}$ is given through a sequence
of functions $B_k:\ \R\longrightarrow \R$, $k\in \N^+$, which are usually
Laurent series in some variable $x$ (see for example \cite{KreimerB}). In what follows, we simply
write the sequence as $\{B_k\}_{k\in \N^+}$.
\begin{definition}\label{bare}
The bare Green function $B\in {\bf B}$, $B:{\cal F}\longrightarrow \R$, associated with 
the sequence of functions $\{B_k\}_{k\in {\bf N_+}}$ is defined recursively
as
$$B(t)=B_{\vert t\vert}\prod_{i=1}^k B(t_i),$$
where $t_1,\cdots,t_k$ are the $d(r)$ subtrees of $t$ rooted
at ${\rm ch}(r)$, and where $\vert t\vert$ denotes the number
of nodes of $t$.
\end{definition}
Notice that the value of $B$ at $t\in {\cal F}$ does not
depend on the linear orders and is independent of the labellings.
When dealing with rooted trees, we will adopt the notation
$t=B_+(t_1,\cdots,t_k)$ for the operation of grafting the rooted
trees $t_1,\cdots,t_k$, that is by considering the tree $t$ obtained
by the creation of a new node $r$ (the root) and then joining
the roots of $t_1,\cdots, t_k$ to $r$.
Bare Green functions appeared also in the probabilistic literature
in specific situations. The basic example, in
algorithmics \cite{Dob,Mah},
 in
numerical analysis (see \cite{ButcherB,ButcherA,Hairer}),
in stochastics
 \cite{Fill,MeirB}
and in physics
 (see for example
\cite{Brouder}) is the {\it tree factorial}, defined by
\begin{definition}\label{Factorial}
Let $t\in {\cal R}$ with 
$t=B_+(t_1,\cdots, t_k)$. Then the
tree factorial $t$! is the functional
$B\in {\bf B}$ defined by
$t!= \vert t\vert \prod_{i=1}^k t_i!$,
associated with the sequence $\{B_k\}$ given by
$B_k\equiv k$.
\end{definition}
\begin{remark}\label{Add1}
It should be pointed out that
the functional acting on trees,
given as $s(t)=\ln(B(t))$, for $B\in {\bf B}$
with $B_k>0$, $\forall k\ge 1$, is
an inductive map or an additive tree functional,
as defined in \cite{berg}. Interestingly, 
$B(t)=1/t!$ is used in \cite{Mah} to define a probability
measure on random search binary trees, and \cite{Dob,Fill} provide
precise asymptotics for
$\ln(t!)$. 
\end{remark}

\section{Generating functions}\label{generating}
We first give some basic results on tree factorials, symmetry factors,
and generating functions associated with bare Green functions.

\begin{definition}\label{Factors}
Let $t\in{\cal R}$. Then
$\alpha(t)$ is the number of rooted labelled trees
$t'\in {\cal L}$ of shape $t\in {\cal R}$, where the shape
of a labelled tree $(r,V,E,l)$ is $(r,V,E)$,
$\kappa(t)$ is the number of rooted plane trees
of shape $t$, and $\sigma(t)$ is the symmetry factor 
of the tree, to be defined later. Moreover, let
$\omega_{\cal L}$ be the weight function associated
with elements of ${\cal L}$, with weights given by
$\psi_k\equiv 1/k!$.
\end{definition}
Notice that $\alpha (t)$ is the
Connes-Moscovici weight in quantum field theory (see 
\cite{Broadhurst,Connes}). The symmetry factor satisfies the recursive
definition:
$$\sigma(\{r\})=1,$$
$$\sigma(B_+(t_1^{n_1},\cdots ,t_k^{n_k}))=n_1!\sigma(t_1)^{n_1}\cdots n_k!\sigma(t_k)^{n_k},$$
where the indices $n_i$ means that $t$ is obtained by grafting $n_1$ times the tree $t_1$, and so on,
where we assume that the $t_i$ are all different as rooted trees.
\begin{lemma}\label{Identities}
Let $t\in {\cal R}$. Then
\begin{equation}\label{symA}
\alpha(t)\sigma(t)=\frac{\vert t\vert !}{t!},
\end{equation}
and
\begin{equation}\label{symB}
\alpha(t) t!=\vert t\vert ! \omega_{\cal L}(t)\kappa(t).
\end{equation}
\end{lemma}
{\it Proof:}
(\ref{symA}) is well known (see for example \cite{ButcherA}).
Suppose that $t\in {\cal R}$ is such that
$t=B_+(t_1^{n_1},\cdots,t_k^{n_k})$, the grafting of 
$n_1$ times the tree $t_1$, and so on,
where we set that the trees $t_1,\cdots,t_k$ are different as rooted trees.
Then
$$\kappa(t)=\frac{(n_1+\cdots+n_k)!}{n_1!\cdots n_k!}\kappa(t_1)^{n_1}\cdots \kappa(t_k)^{n_k}.$$
Using the recursive definition of $\omega(t)$ and the definition of
$\oml$, we have
$$\oml (t)=\frac{1}{(n_1+\cdots+n_k)!}\oml (t_1)^{n_1}\cdots \oml (t_k)^{n_k}.$$
Therefore
$$\frac{1}{\oml (t)\kappa(t)}=n_1!\cdots n_k! (\frac{1}{\oml(t_1)\kappa(t_1)})^{n_1}\cdots (\frac{1}{\oml(t_k)\kappa(t_k)})^{n_k},$$
and the results follows from
 the recursive definition of the symmetry factor.

\endproof
Then
\begin{equation}\label{use}
\sum_{t\in {\cal F}_n}B(t)\omega (t)=\sum_{t\in {\cal R}_n}B(t)\frac{\omega}{\oml}(t)\alpha(t)\frac{t!}{\vert t\vert!}
\end{equation}
where we have used (\ref{symB}) of Lemma \ref{Identities}.

\noindent Consider the generating function
\begin{equation}\label{Generating}
Y(z)=\sum_{n\in \N^+ }
\frac{z^n}{n!}\sum_{t\in {\cal R}_n}\alpha(t)B(t)t!\omega(t)/\oml(t).
\end{equation}
 Given $t\in {\cal R}$, the ratio $\omega/\oml$ is associated with the weight sequence
$\bar\psi_k\equiv \psi_k k!$; using the expansion $\psi(z)=1+\sum_{k\ge 1}\psi_k z^k=1+\sum_{k\ge 1}(\bar\psi_k /k!)z^k$,
we see that $\bar\psi_k\equiv \psi^{(k)}(0)$. Consider the {\it elementary differentials} $\delta$  (see Section \ref{Runge})
  defined by
\begin{definition}\label{elementary}
$$\delta_{\{*\}}=1,\ \delta_t=\psi^{(k)}(0)\prod_{i=1}^k \delta_{t_i},\ \ \frac{\omega}{\oml}=\delta,$$
when $t=B_+(t_1,\cdots, t_l)$, where * denotes the tree of a single node.
For a map $a:\ {\cal R}\cup \{\emptyset\}\longrightarrow \R$,
    a formal power series of the form
$Y(z)=a(\emptyset)y_0+\sum_{t\in {\cal R}}z^{\vert t\vert}a(t)\delta_t \alpha (t) /\vert t\vert!$ is called a B-series
 \cite{Wanner,Hairer}.
\end{definition}
\begin{remark}\label{Add2}
When $B(t)=t!$,  the series $Y$ is given by
  $$Y(z)=\sum_{t\in {\cal L}}(\omega(t)/\oml(t))z^{\vert t\vert}/\vert t\vert !.$$
 Set
  $\phi_k=\psi_k k!$, $\forall k$, and consider the degree function $\phi(z)=1+\sum_{k\ge 1}(\phi_k/k!)z^k$.
  Following \cite{berg}, $Y$ solves $Y'=\phi(Y)$ (see also \cite{MeirB}). We shall see in the next section
  that it is a natural consequence of B-series expansions of solutions of ordinary
  differential equations.
  \end{remark}


\section{Runge-Kutta methods for functionals over trees}\label{Runge}
   Consider a dynamical system on $\R$
  $$\ds X(s)=F(X(s)),\ \ X(s_0)=X_0,$$
  for some smooth $F:\R\longrightarrow \R$.
   The solution of this dynamical system
  has a B-series expansion of the form
  $$X(s)=X_0+\sum_{t\in {\cal R}}\frac{(s-s_0)^{\vert t\vert}}{\vert t\vert!}\alpha(t)\delta_t(s_0),$$
  where the elementary differentials $\delta$ is defined recursively by
  $$\delta_{\{r\}}=f(s_0),\
  \delta_t=
  \frac{\partial^k F }{\partial^k s}\delta_{t_1}\cdots\delta_{t_k},$$
  when $t=B_+(t_1,\cdots, t_k)$. These kinds of expansions have been treated in great detail in \cite{ButcherB}
  and \cite{ButcherA}
  and developped independently in combinatorics (see for example \cite{LerouxA,LerouxB}).
 Suppose that $s_0=0$ for simplicity.
   Butcher considered
  what happens with numerical approximations of the exact solution,
  the Runge-Kutta methods, which are themselves B-series
   \cite{Wanner,Hairer}; here we focus
  on specific B-series, which are associated to bare Green functions.
  Let $B\in {\bf B}$ be such that
there exists a power series
$$L(z)=\sum_{m\ge 0}L_m z^m,$$
with
\begin{equation}\label{Choice}
B_k =
\frac{L(k)}{k}, \ \forall k\in \N^+.
\end{equation}
Bare Green functions satisfying (\ref{Choice})
are used in practical situations in quantum field theory
(see  \cite{Broadhurst}, \S~4 and
\cite{Brouder} , \S~6.1).
Consider Euler's operator
$\theta=z ({\rm d}/{\rm d}z)$,
with
$P(\theta)(z^n)=P(n)z^n,\ \forall n\in \N$,
for any polynomial $P$, and consider the formal
operator $L(\theta+1)$ acting
on monomials as
\begin{eqnarray*}
L(\theta+1)(z^n)&=&
\sum_{m\ge 0}L_m (\theta+1)^m (z^n)=
\sum_{m\ge 0}L_m (n+1)^m z^n\\
&=& L(n+1)z^n.
\end{eqnarray*}
Given a power series $Y(z)=\sum_{m\ge 0}a_m z^m$ converging
for $\vert z\vert\le 1$, we can define
$L(\theta+1)(Y)(z):=\sum_{m\ge 0}a_m L(m+1)z^m$, which
converges for $\vert z\vert\le 1$ when
the sequence $(L(k))_{k\ge 1}$ grows subexponentially. We will not
focus on convergence questions here, and work at the formal level.
Let $B$ be a bare Green function with weights
$(B_k)_{k\ge 1}$, such that
(\ref{Choice}) holds for some power series $L$.
It should be pointed out that \cite{Broadhurst,Brouder,Kreimer} deal with
the master function $L$, but  do not
give explicitely an equation for $Y$. The next Theorem
provides an equation; its proof uses explicitely B-series.
\begin{theorem}\label{DifferentialSys}
The formal power series
\begin{equation}\label{Bseries}
Y(z)=\sum_{t\in {\cal R}}\frac{z^{\vert t\vert}}{\vert t\vert!}
\alpha(t) t! B(t)\delta_t,
\end{equation}
solves   $Y'=L(1+\theta)\psi(Y)$.
\end{theorem}
 {\it Proof:}

  $$\psi(Y(z))=\sum_{k\ge 0}\frac{\psi^{(k)}(0)}{k!}\sum_{(t_1\cdots t_k)\in {\cal R}^k}
  \frac{z^{\sum_i \vert t_i\vert}}{\vert t_1\vert!\cdots \vert t_k\vert!}     \prod_{i=1}^k
     \alpha(t_i)
  B(t_i)t_i!
  \delta_{t_i}.$$
   For given $(t_1\cdots t_k)\in {\cal R}^k$, set
   $t=B_+(t_1,\cdots, t_k)$. Then
   $\sum_i \vert t_i\vert=\vert t\vert -1$, $\psi^{(k)}(0)\delta_{t_1}\cdots\delta_{t_k}=\delta_t$,
   and $B(t_1)\cdots B(t_k)=B(t)/B_{\vert t\vert}$.
 The associated term becomes
   $$  z^{\vert t\vert -1}\frac{B(t)}{B_{\vert t\vert}}\delta_t \frac{\alpha(t_1)\cdots\alpha(t_k)}{\vert t_1\vert!\cdots \vert t_k\vert!}
   \frac{t!}{\vert t\vert}.$$
 Next, every rooted tree $t\in {\cal R}$ can be decomposed uniquely as 
 $t=B_+$ $(t_1^{n_1},\cdots, t_m^{n_m})$,
   meaning that $t$ is obtained by grafting $n_1$ times $t_1$ and so on, where the $t_i$ are different as rooted trees,
   with $k=n_1+\cdots+n_m$.
   Collecting the terms associated with $t$, we get the contribution
   $$
    \frac{z^{\vert t\vert -1}}{k!}\frac{B(t)}{B_{\vert t\vert}}\delta_t\frac{t!}{\vert t\vert}
   \sum_{(t'_1\cdots t'_k)\in {\cal R}^k}^*
       \frac{\alpha(t' _1)\cdots\alpha(t' _k)}{\vert t'_1\vert!\cdots\vert t'_k\vert!},$$
       where * means that the sum is taken over all the collections
       $(t'_1\cdots t'_k)\in {\cal R}^k$ such that $t=B_+(t'_1,\cdots ,t'_k)$.
       The above sum reduced then to
       $$\frac{(n_1+\cdots+n_m)!}{n_1!\cdots n_m!}\frac{\alpha(t_1)^{n_1}\cdots \alpha(t_m)^{n_m}}{k! \vert t_1
       \vert !^{n_1}\cdots \vert t_m\vert  !^{n_m}}=
       \frac{1}{n_1!}(\frac{1}{\sigma(t_1)t_1!})^{n_1}\cdots\frac{1}{n_m!}(\frac{1}{\sigma(t_m)t_m!})^{n_m},$$
       where we use the first identity of Lemma \ref{Identities}. Using the recursive definition of the symmetry factor $\sigma$, we obtain
       \begin{eqnarray*}
       \sum_{(t'_1\cdots t'_k)\in {\cal R}^k}^*
       \frac{\alpha(t' _1)\cdots\alpha(t' _k)}{k!\vert t'_1\vert!\cdots\vert t'_k\vert!}&=&
       \frac{1}{t_1!^{n_1}\cdots t_m!^{n_m}\sigma(B_+(t_1,^{n_1}\cdots, t_m^{n_m}))}\\
       &=&\frac{\vert t\vert}{t!}\frac{1}{\sigma(t)}=\frac{\vert t\vert \alpha(t)}{\vert t\vert!}.
       \end{eqnarray*}
       We thus get that the contribution associated with $t\in  {\cal R}$ is given by
       $$
          \frac{z^{\vert t\vert -1}}{k!}\frac{B(t)}{B_{\vert t\vert}}\delta_t\frac{t!}{\vert t\vert}\frac{\vert t\vert \alpha(t)}{\vert t\vert!}
                 = \frac{ z^{\vert t\vert -1}}{\vert t\vert!} \frac{B(t)}{B_{\vert t\vert}}\alpha(t)\delta_t t! .$$
      Therefore
      \begin{eqnarray*}
      L(\theta+1)\psi(Y)&=&\sum_{t\in {\cal R}}\frac{B(t)}{B_{\vert t\vert}}\alpha(t)t!\delta_t \frac{L(\theta+1)(z^{\vert t\vert-1})}{\vert t\vert!}\\
                   &=&\sum_{t\in {\cal R}} \frac{B(t)}{B_{\vert t\vert}}\alpha(t)t!\delta_t    \frac{L(\vert t\vert)z^{\vert t\vert-1}}
                   {\vert t\vert!}\\
                   &=&    \sum_{t\in {\cal R}} \frac{B(t)}{B_{\vert t\vert}}\alpha(t)t!\delta_t    \frac{B_{\vert t\vert}\vert t\vert
                   z^{\vert t\vert-1}}{\vert t\vert!}\\
                   &=&\sum_{t\in {\cal R}}\frac{z^{\vert t\vert -1}}{(\vert t\vert -1)!}B(t)t!\delta_t=\frac{{\rm d}Y}{{\rm d}z}.
     \end{eqnarray*}

 \endproof
\begin{remark}\label{Add3}
As we have observed in Remark \ref{Add1},
$s(t)=\ln(B(t))$ is an inductive map when
the weights $B_k$ are positive. It turns out
that the exponential generating function associated
with $s$ can be given as an integral transform
for varieties of increasing trees (see for example Section \ref{s.intro}). This is the topic
of \cite{berg}.
\end{remark}
\begin{example}\label{GenCat}
\end{example}
\noindent When $L(z)=z$, with $B_k\equiv 1$, and
$\psi(z)=1/(1-z)$, one has
$\sum_{t\in {\cal F}_{n+1}}B(t)=C_n$, the 
Catalan number of order $n$, with
$C_n={2n\choose n}/(n+1)$. Then
$Y(z)=z\sum_{n\ge 0}z^n C_n$ is solution
of the differential equation
$Y'(z)=L(1+\theta)(1/(1-Y(z)))$, that is
$Y'(z)=(z/(1-Y(z)))'$. The unique solution
with $Y(0)=0$ satisfies $Y(z)=z/(1-Y(z))$,
or $Y(z)=(1-\sqrt{1-4z})/2$, corresponding
to a well known result.

\endproof
\noindent \cite{Brouder} , \S~5.3, considers the case
 where $B(t)=(1/t!)^2$, which is not of the form
 given in (\ref{Choice}): in this situation,
 $B_k=1/k^2$, with $L(z)=1/z$. The solution
    is obtained by using
   the stucture
 of the so-called Butcher's group of
 B-series  (that is series of the form (\ref{Bseries}),
 where the group structure in given in
 \cite{Wanner,Hairer})
   by tensoring
 known B-series:
    \begin{example}\label{m2}
    \end{example}
\noindent Consider the bare functional given by
$B_k\equiv 1/k^2$, with $B(t)=1/t!^2$.
Following Brouder, the associated B-series, as given in (\ref{Bseries}), is solution of the second
order differential equation
$$z Y'' +Y'=\psi(Y).$$
When $\psi(z)=\exp(z)$, the
 solution is given by
$$Y(z)=-2\ln(1-z/2)=\sum_{n=1}^{\infty}\frac{z^n}{n}\frac{1}{2^{n-1}},$$
giving
$$\sum_{t\in {\cal R}_n}\frac{\alpha(t)}{t!}=\frac{(n-1)!}{2^{n-1}}.$$

\endproof
 We study   the
 general moment problem $B(t)=(1/t!)^{ l +1} $, $l\in \N$,
 by working directly on a suitable differential equation
 as follows: the operator
 $L(\theta+1)$ takes the form
 $L(\theta+1)=1/(\theta+1)^l$.
 Assume that the differential operator $L(\theta+1)$ is invertible.
 Then the formal
  systems
 becomes
 \begin{equation}\label{EquiSysA}
 L(\theta+1)^{-1}\dz Y =\psi(Y).
 \end{equation}
 Consider again the second moment problem for the tree factorial,
 with $B_k\equiv 1/k^2$ and $L(k)\equiv 1/k$. Choose $L$
 such that $L(z)=1/z$; the inverse operator might be equal to
 $L(\theta+1)^{-1}=\theta+1$, and, if this is the case,
 $$(\theta+1)\dz Y=\psi(Y),$$
 with
 $(\theta+1)(d/dz)
 =z({\rm d^2}/{\rm dz^2})+({\rm d}/{\rm d}z)$,
 see Example
 \ref{m2}.

 More generally, if one considers the moment of order $l+1
 \in \N$ of the inverse
 tree factorial, the choice $L(k)=1/k^l$
  should give
 $(\theta +1)^{l}\dz Y(z)=\psi(Y)$.
Our result, Theorem \ref{Inversion} below
  shows that the formalism of inversion is correct
 in term of power series. This result sheds some light
 and extends the computations done in \cite{Brouder}
 for the second moment, and its proof avoids computations
 in the Butcher's group.
 
\begin{theorem}\label{Inversion}
The B-series $Y(z)$ associated with
the moment of order $(l+1)$ of the inverse tree factorial satisfies
the differential equation
\begin{equation}\label{InverseEquation}
  (\theta+1)^l \dz Y=\psi(Y).
  \end{equation}
  \end{theorem}
  {\it Proof:}
  Let
  $$Y(z)=\sum_{t\in {\cal R}}\frac{z^{\vert t\vert}}{\vert t\vert!}\alpha(t)\frac{1}{t!^{l+1}}t!\delta_t.$$
  Then
  $$\psi(Y(z))=\sum_{k\ge 0}\frac{\psi^{(k)}(0)}{k!}\sum_{(t_1\cdots t_k)\in {\cal R}^k}
  \frac{z^{\sum_i \vert t_i\vert}}{\vert t_1\vert!\cdots \vert t_k\vert!}\frac{\alpha(t_1)\cdots\alpha(t_k)}{(t_1!\cdots t_k!)^l}\delta_{t_1}\cdots\delta_{t_k}.$$
   For given $(t_1\cdots t_k)\in {\cal R}^k$, set
   $t=B_+(t_1,\cdots, t_k)$. Then
   $\sum_i \vert t_i\vert=\vert t\vert -1$, $\psi^{(k)}(0)\delta_{t_1}\cdots\delta_{t_k}=\delta_t$,
   and $(t_1!\cdots t_k!)^l=t!^l/\vert t\vert^l$. The associated term becomes
   $$\frac{z^{\vert t\vert-1}\vert t\vert^l}{t!^l}\delta_t\frac{\alpha(t_1)\cdots\alpha(t_k)}{\vert t_1\vert!\cdots\vert t_k\vert!}.$$
   Proceeding as in the proof of Theorem \ref{DifferentialSys},
      we
       get that the contribution associated with $t\in  {\cal R}$ is given by
       $$  \frac{z^{\vert t\vert-1}\vert t\vert^l}{t!^l}\delta_t  \frac{\vert t\vert \alpha(t)}{\vert t\vert!}
       =\frac{z^{\vert t\vert-1}}{ (\vert t\vert-1)!}\frac{\vert t\vert^l}{t!^l}\alpha(t)\delta_t.$$
       On the other hand,
       \begin{eqnarray*}
       (\theta+1)^l\dz Y(z)&=&(\theta+1)^l\sum_{t\in {\cal R}}\frac{z^{\vert t\vert -1}}{(\vert t\vert -1)!}
       \alpha(t)\frac{1}{t!^l}\delta_t\\
       &=&\sum_{t\in {\cal R}}\frac{\vert t\vert^l z^{\vert t\vert -1}}{(\vert t\vert-1)!}\alpha(t)\frac{1}{t!^l}\delta_t.
       \end{eqnarray*}

       \endproof

In the next section, we show that traces of certain
products of Wigner matrices (see for example \cite{Sos}) provide
natural examples of bare Green functions.
\section{Wigner processes\label{WignerProcesses}}

\begin{definition}
\label{d.goe}
The $N$-dimensional random matrices
$
\Gamma_N:=(\gamma_{i,j})_{1\le i,j\le N}
$
are called Wigner matrices
of variance $\beta^2$ if the following holds.
\begin{itemize}
\item
Each $\Gamma_N$ is symmetric, that is,
$\gamma_{i,j}=\gamma_{j,i}$.
\item
For $i\le j$, the random variables $\gamma_{i,j}$ are independent and
centered.
\item
For $i\neq j$, $\E(\gamma_{i,j}^2)=\beta^2$.
\item
For any $k\ge 2$, $\E(|\gamma_{i,j}|^k)\le c_k$, where
$c_k$ is independent of $i\le j$.
\end{itemize}
\end{definition}

\begin{definition}
\label{d.goep}
The sequence $\Gamma_N(k):=(\gamma_{i,j}(k))_{1\le i,j\le N}$
of $N$-dimensional random matrices, indexed by $k\ge 1$,
is called a
Wigner process of variance $\beta^2$ and
correlation function
 $r$,
 $r(k,k)=1$,
 $\vert r(k,m)\vert\le 1$ and
 $r(k,m)=r(m,k)$
 if the following holds.
\begin{itemize}
\item
Each $\Gamma_N(k)$ is a Wigner
matrix of variance $\beta^2$ in the sense of definition~\ref{d.goe}.
\item
For $i\le j$, each process $(\gamma_{i,j}(k))_k$ is
independent of the others.
\item
For $i\neq j$, the process $(\gamma_{i,j}(k))_k$ is $r$-correlated,
that is,
for any $k\ge m$,
\begin{equation}
\label{e.corr}
\E(\gamma_{i,j}(k)\gamma_{i,j}(m)):=\beta^2\,r(k,m).
\end{equation}
\end{itemize}
A Wigner process is stationary when $r$ is such that
$r(k,m)=r(\vert k-m\vert)$.
\end{definition}
Let $D_N$ be a sequence of random diagonal matrices,
with independent and identically distributed entries  of law $\mu$, having finite moment
$\mu_k=\mu(X^k)$, $k\ge 1$, with
$\mu_1=1$.
Let
$$
Q_{k}^N:=N^{-k/2}D_N\prod_{m=1}^{k}(\Gamma_N(m)D_N),
$$
 and set
$$ B_k^N(r)=
N^{-1}\,\E(\trr(Q_{k}^N)).
$$
\bigskip

\noindent{\bf Involutions, Dyck paths and rooted plane trees}
\bigskip

For $k\ge 1$, $[k]:=\{1,2,\ldots,k\}$,
$\I(k)$ is the set of the involutions of $[k]$ with no fixed point,
$\J(k)$ is the subset of $\I(k)$ of the involutions $\sigma$ with no crossing.
This means that the configurations
$
i<j<\sigma(i)<\sigma(j)
$
do not appear in $\sigma\in\J(k)$.
Let $i\in\crr(\sigma)$ denote the fact that $i<\sigma(i)$.
Let $\D(2k)$ be the set of the Dyck paths of length $2k$, that is, of
the sequences 
$c:=(c_n)_{0\le n\le 2k}$ of
nonnegative integers such that
$
c_0=c_{2k}=0$,
$
c_{n}-c_{n-1}=\pm 1$,
$ n\in[2k]$.
Thus, exactly $k$ indices $n\in[2k]$ correspond to ascending steps
$(c_{n-1},c_n)$, that is, to steps when $c_n=c_{n-1}+1$.
We denote this by $n\in\asc(c)$.
The $k$ others indices correspond to descending steps,
that is, to steps when $c_n=c_{n-1}-1$,
and we denote this by $n\in\desc(c)$.
We make use of bijections between $\D(2k)$ and $\J(2k)$ \cite{Biane}.
If $c\in\D(2k)$,
$\phi(c):=\sigma\in\J(2k)$ is an involution which maps each element of
$\desc(c)$ to a smaller element of $\asc(c)$.
Thus, $\crr(\sigma)=\asc(c)$.
More specifically,
if $n\in\desc(c)$, $\sigma(n)$ is the greatest $m\le n$ such
that $(c_{m-1},c_m)=(c_n,c_{n-1})$.
Finally, the set ${\cal D}(2k)$ is in bijection with
${\cal F}_{k+1}$, the set of rooted plane trees on $k+1$ nodes,
where the bijection is given by the walk on the
tree from the right to the left (see for example \cite{Takacs}). Let $\sigma_t$ denote the
involution of ${\cal J}(2(\vert t\vert-1))$ corresponding
to $t\in {\cal F}$.  Given $t\in {\cal F}_{k+1}$, consider
the walk on $t$ from the right to the left: every
edge $(v,w)$ with $w\in {\rm ch}(v)$, is crossed at some instant
$s_v\in [2k]$ as $(v\to w)$ and at a later time
$s_w\in [2k]$ as $(w\to v)$. Clearly,
$s_w=s_v+2(\vert t_w\vert -1)+1$, where
$t_w$ is the subtree of $t$ rooted at node $w$,
that is the subgraph of $t$ induced by the nodes
$u$ with $u\ge w$.
$\sigma_t$ is such that
$\sigma_t(s_v)=s_w$ and vice versa.

\begin{center}
\includegraphics{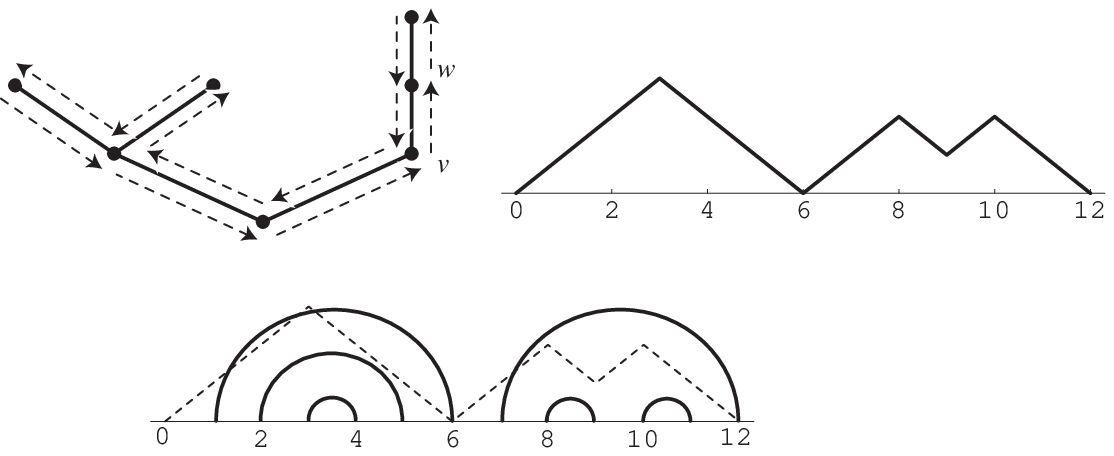}
\end{center}

\noindent Fig 1. Bijections between ${\cal F}_{k+1}$, ${\cal D}(2k)$ and $\J
(2k)$ ($s_v=2$ and $s_w=5$).

\begin{proposition}\label{SpecialBare}
Assume that the covariance $r$ is such that
$r(l,m)=r(\vert l-m\vert)$. Then,
the functional $B^r\in {\bf B}$
given by
the weights
$$B^r_k=\beta^2 r(2k-1),\ \forall k\ge 1,$$
is such that
$$\frac{B^r(t)}{B^r_{\vert t\vert}}
=\prod_{i\in cr(\sigma_t)}( \beta^2
r(i,\sigma_t(i))).$$
\end{proposition}
{\it Proof:}  Let $t\in {\cal F}$. Let $s_v<s_w$ be the
instants where the oriented edges
$(v\to w)$ and $(w\to v)$, $w\in {\rm ch}(v)$, are
crossed during the walk on the tree.
$r(s_v,\sigma_t(s_v))=r(s_w-s_v)=r(2\vert t_w\vert-1)$,
and thus $\beta^2 r(s_v,\sigma_t(s_v))=B^r_{\vert t_w\vert}$.
Finally, $\prod_{i\in cr(\sigma_t)}\beta^2 r(i,\sigma_t(i))
=\prod_{w\ne {\rm root}}B^r_{\vert t_w\vert}=B^r(t)/B^r_{\vert t\vert}$,
 as required.

 \endproof
\medskip

As we have just seen, every
Wigner process
with covariance $r$ such that
$r(l,m)=r(\vert l-m\vert)$ produces a
bare Green function $B^r\in {\bf B}$. The
converse is not true, that is, there exists $B\in {\bf B}$
such that $B$ is not of the form $B=B^r$ for some
covariance function $r$.    Set
${\bf B}^w=\{B\in {\bf B};\ \exists \hbox{ a covariance }r\hbox{ with }B=B^r\}$.

Let $\psi_\mu$ be the generating function of the
weight sequence $\psi_k=\mu_{k+1}$, and let
$\omega_\mu(t)$, $t\in {\cal F}$ be the associated
weight function.
\begin{theorem}
\label{t.wig}
Let $(\Gamma_N(k))_{k\ge 1}$ be a stationary Wigner
process of covariance function $r$ and variance $\beta^2$, and let
$D_N$ be a sequence of random diagonal matrices,
independent of the Wigner process,
with i.i.d. entries $\lambda_j$ of law $\mu$,
with $\mu_1=\mu(\lambda)=1$
and  finite moments $\mu_k=\mu(\lambda ^k)$,
$\forall k$.
 Then
$$B_{2k}^N(r)\longrightarrow B_{2k}(r)=
\frac{1}{B_{k+1}^r}\sum_{t\in {\cal F}_{k+1}}B^r(t)\omega_\mu(t),$$
and $B_{2k+1}^N(r)\longrightarrow 0$, $N\to\infty$.
Assume that the covariance is such that there exists a power series $L^r(z)$ with
$B^r_k= L^r(k)/k$, $\forall k$. Then the formal power series
$$Y(z)=\sum_{k\ge 1}z^k  B_k^r B_{2(k-1)}(r),$$
solves
$$Y'=L^r(\theta +1)\psi_\mu(Y).$$
Moreover
\begin{equation}\label{FondLemma}
\sum_{k\ge 1}z^k B_{2(k-1)}(r)=z\psi_\mu(Y).
\end{equation}
\end{theorem}
  \begin{example}\label{TreeFAct}
  \end{example}
\noindent Let $B(t)=1/t!$. If a tree   $t$
 has $n$ nodes
 and $n-1$ edges,    then the requirement
 $B_n=1/n$ is satisfied iff $\beta^2
 r(2n-1)=1/n$, that is
 $r$ must be such that
 $ \beta^2
 r(k)=2/(k+1)$, $k\in 2\N+1$.
 By construction, $r(0)=1$ and therefore
 $\beta^2=2  $. $1/(x+1)$ is positive
 definite, which implies that
 $B(t)=1/t!$ is element of ${\bf B}^w$.
 Next, from Theorem \ref{Inversion}, the generating function $Y(z)=\sum_{t\in {\cal F}}
 z^{\vert t\vert}B^r(t)\omega_\mu(t)$ is solution of the system
 $({\rm d}/{\rm d}z)Y(z)=\psi_\mu(Y)$. Assume that
 $\mu$ is the point mass $\delta_1$, that is each matrix
 $D_N$ is the identity matrix of size $N$,
 with $\psi_\mu(z)=1/(1-z)$. The solution
 of the system is $Y(z)=1-\sqrt{1-2z}=2\tilde Y(z/2)$,
 where $\tilde Y$ is the series given in Example
 \ref{GenCat}.
On the other hand,
 Proposition \ref{SpecialBare} and Theorem \ref{t.wig}
 show that $Y(z)=\sum_{k\ge 1}z^k B_k^r B_{2(k-1)}(r)$. Therefore
 the limiting mean normalized trace $B_{2k}(r)$ of the product of correlated
 random matrices $\prod_{m=1}^{2k}\Gamma_N(m)$ is such that
 $B_{2k}(r)={\rm E}(Z^{2k})/k!$, where $Z$ denotes a
 normal N(0,1) random variable.

 \endproof

\begin{example}\label{GenCat2}
\end{example}
\noindent Consider as in Example \ref{GenCat} the special case
where $L(z)=z$. The associated inductive parameter
(see Remark \ref{Add3}) is the tree size.
The covariance $r$ is constant with
$r(k)\equiv 1$, and $B^r_k\equiv 1$. Then the generating
function $Y$ is solution of the fixed point equation
$Y(z)=z\psi_\mu(Y(z))$
(either by Theorem \ref{DifferentialSys} or by
(\ref{FondLemma})). Notice that in this situation,
$\Gamma_N(m)\equiv \Gamma_N(1)$, and thus
$B_k^N(r)$ describes the mean normalized moment
 of the spectral measure of the random matrix $D_N(\Gamma_N(1)D_N)^k$.
This example can be
extended by considering
$L(z)=z\rho^z$,
for some
$0<\rho\le 1$. When 
$D_N$ is the identity matrix, $Y(z)$ is related
to the Rogers-Ramanujan continued fraction
 \cite{Mazza}, and corresponds to the generating
function associated with  path length, see \cite{berg,Takacs}.

\endproof
\bigskip

\noindent {\it Proof of Theorem \ref{t.wig}:} The first part is
a generalization of
Theorem 1 of \cite{Mazza}. Set $\tilde\gamma_{ij}(m)=\gamma_{ij}(m)\lambda_j$,
and $\tilde\Gamma_N(m)=\Gamma_N(m)D_N$. The mean normalized trace
adds the contributions $E(i)={\rm E}(\lambda_{i_0}\tilde\gamma_{i_0 i_1}\cdots \tilde\gamma_{i_{k-1}i_k})$,
for paths $i=(i_l)_{0\le l\le k}$, with $i_l\in [N]$ and
$i_0=i_k$. The $\tilde\gamma_{ij}$ are centered, so that
any edge $(i,j)$ appearing once appears at least twice.
Given $i$, define $\varepsilon_1=1$ and
$\varepsilon_l = +1$ when $i_l\not\in\{i_0,\cdots,i_{l-1}\}$,
and $\varepsilon_l=-1$ otherwise, and consider the walk
$c=(c_l)$ defined by $c_l=\sum_{j=1}^l \varepsilon_j$,
with $c_k\le 0$. The support of $i$ is
$s(i)=\{i_l;\ 0\le l\le k\}$, of
size $s=\vert s(i)\vert$, with $s\le 1+k/2$. The contribution
$E(i)$ is independent of the labels $i_l$; they are
$N(N-1)\cdots (N-s+1)$ labellings giving the same walk
$c$, with the same contribution. Thus, the normalization
$N^{-(1+k/2)}$ shows that the only walks surviving in the large
$N$ limit are those with $s=1+k/2$. This shows that 
$B_k^N(r)\to 0$ when $k$ is odd. Concerning $B_{2k}^N(r)$,
$s=1+k$ means that every edge occuring in the path
occurs exactly twice, in opposite directions. $c$
is a Dyck path of ${\cal D}(2k)$; let $t\in {\cal F}$
be the associated rooted plane tree, with involution
$\sigma_t$. Using the right to left walk on $t$
and the independence of the random variables,
 the contribution $E(i)$
of any path leading to $c$ or $t$ is
$E(i)=\prod_{m\in cr(\sigma_t)}{\rm E}(\gamma(m)\gamma(\sigma_t(m)))
{\rm E}(\prod_{v}\lambda_v^{d(v)+1})$
where $d(v)=\vert {\rm ch}(v)\vert$. From Proposition
\ref{SpecialBare}, one obtains 
$E(i)=(B^r(t)/B_{k+1}^r)\prod_v \mu_{d(v)+1}$, with
$B_{2k}(r)=\sum_{t\in {\cal F}_{k+1}}(B^r(t)/B_{k+1}^r)\prod_v \mu_{d(v)+1}$,
as required. (\ref{FondLemma}) is a consequence of the multiplicative
form of bare Green functions and of Lemma 1.9, chap. III.1 of \cite{Hairer}.

\endproof

 \medskip

These results show that the elements of ${\bf B}^w$ appear naturally 
in the computation of 
 normalized traces of products of large random matrices
 (see for example \cite{Speicher}). In the next Section we illustrate
 B-series  by considering
triangular operators from free probability.


\section{On Dykema-Haagerup triangular operator\label{Triangular}}

Let ${\cal B}$ be an algebra and ${\cal A}$ be
a ${\cal B}$ bi-module. Let $\kappa:\ \A\ {\rm x}\ \A\longrightarrow \B$
be a bilinear map. We follow \cite{Sni} by defining the product
$a_1\bu a_2=\kappa(a_1,a_2)$, $a_1,a_2\in\A$, and setting
$$i)\ (ba_1)\bu a_2=b(a_1\bu a_2),$$
$$ii)\ (a_1b)\bu a_2=a_1\bu(ba_2),$$
$$iii)\ a_1\bu (a_2b)=(a_1\bu a_2)b.$$
Let $\sigma\in {\cal J}(2n)$ be an involution
of $[2n]$ without fixed point and without crossing.
Given a word $a=a_1\cdots a_{2n}$ in $\A$, $\sigma$
induces parentheses on $a$, and the preceedings rules
permit the evaluation of this parenthized word. This
extends to a map $\kappa_\sigma$ on $\A^{2n}$.
Sniady defines such maps to prove a conjecture
of Dykema and Haagerup on generalized circular elements.
Let $(\B\subset A,{\rm E})$ be
an operator valued probability space, that is
$\A$ is a unital *-algebra, $\B\subset A$ an unital
*-subalgebra and ${\rm E}:\A\longrightarrow \B$ be
a conditional expectation
(linear, ${\rm E}(1)=1$, and
${\rm E}(b_1 ab_2)=b_1{\rm E}(a)b_2$,
$\forall b_1, b_2\in \B,\ a\in\A$).
\begin{definition}\label{Generalized}
$T\in \A$ is a generalized circular element
if there is a bilinear map
$\kappa$ satisfying the rules i), ii) and iii)
such that
$${\rm E}(b_1T^{s_1}b_2T^{s_2}\cdots b_{2n}T^{s_{2n}})=
\sum_{\sigma\in {\cal J}(2n)}\kappa_\sigma(b_1T^{s_1},\cdots,b_{2n}T^{s_{2n}}),$$
$${\rm E}(b_1T^{s_1}b_2T^{s_2}\cdots b_{2n+1}T^{s_{2n+1}})=0,$$
$\forall b_1,\cdots,b_{2n+1}\in\B$ and 
$\forall s_1,\cdots,s_{2n+1}\in \{1,*\}$.
\end{definition}
The triangular operator $T$ of Dykema and Haagerup
is obtained from $\B=\C[x]$, the *-algebra of
complex polynomials of one variable by setting
$$[\kappa(T,bT^*)](x)=\int_x^1 b(s){\rm d}s,$$
$$[\kappa(T^*,bT)](x)=\int_0^x b(s){\rm d}s,$$
$$[\kappa(T,bT)](x)=[\kappa(T^*,bT^*)](x)=0.$$
$T$ is the limit for the convegence of
*-moments of large upper triangular random matrices $T_N$
(\cite{Dykema}).
Define a trace $\tau$ as (see \cite{Sni})
$$\tau(a)=\tau({\rm E}(a)),\ \tau(b)=\int_0^1 b(s){\rm d}s.$$

In what follows, we use P-series (where P stands for
partitioned differential systems, see \cite{Wanner}).
We follow \cite{Brouder}, and adapt his notations to
P-series. Given some function $\psi$, and two kernels
$(a^x(u,v))_{u,v\in [0,1]}$ and
$(a^y(u,v))_{u,v\in [0,1]}$, consider the iterated integrals
$\phi^x_u$ and $\phi^y_u$ which are functionals over 
${\cal R}$ defined by $\phi^x_u(*)=\phi^y_u(*)=1$, and,
for $t=B_+(t_1,\cdots,t_k)$,
$$\phi^x_u(t)=\prod_{i=1}^k \int_0^1 a^x(u,v)\phi^y_v(t_i){\rm d}v,$$
$$\phi^y_u(t)=\prod_{i=1}^k \int_0^1 a^y(u,v)\phi^x_v(t_i){\rm d}v.$$
\begin{lemma}\label{Moments}
Let $a^x(u,v)={\rm I}_{[0,u]}(v)$
and $a^y(u,v)={\rm I}_{[u,1]}(v)$. Then
$$
\tau(TT^*)^n=\sum_{t\in {\cal F}_{n+1}}\int_0^1 \phi^x_v(t){\rm d}v=
            \sum_{t\in {\cal F}_{n+1}}\int_0^1 \phi^y_v(t){\rm d}v.
$$
\end{lemma}
{\it Proof:} The word
$W=(TT^*)\cdots (TT^*)$ is of the generic form
with $b_1=\cdots b_{2n}=1$ (Definition
\ref{Generalized}). Let $t\in {\cal F}_{n+1}$
with associated involution $\sigma_t$
(see Section \ref{WignerProcesses}). Let
$s_v$ and $s_w$ be the instants
where the walk on $t$ crosses
the oriented edges $(v\to w)$
and $(w\to v)$, with $w\in {\rm ch}(v)$. We colour
these edges by giving colour '1' to
$(v\to w)$ when the symbol in $W$
located at position $s_v$ is $T$,
and give the colour '*' otherwise. 
Clearly, both edges have different colours, and
the elements of
the set of edges
$\{(v\to w);\ w\in {\rm ch}(v)\}$
(the children of $v$ in $t$)
have the same colour. The result is then
a consequence of the definition of the product
with the rules i), ii) and iii).

\endproof
\begin{remark}
Iterated integrals are naturel objects
to consider in the setting of Butcher's Theory.
For example, in the framework of Theorem
\ref{DifferentialSys}, the iterated integrals
$\phi_u(t)$ defined by
$\phi_u(t)=\prod_{i=1}^k \int_0^u L(\theta+1)(\phi_v(t_i)){\rm d}v$,
when $t=B_+(t_1,\cdots,t_k)$,
are such that $\phi_1(t)=B(t)$, $\forall t\in {\cal F}$.
\end{remark}

\begin{proposition}\label{IntegralSystem}
The P-series
$$X_u(s)=X_0+\sum_{t\in {\cal R}}\frac{s^{\vert t\vert}}{\vert t\vert!}
\alpha(t)t!\delta_t \int_0^1 a^x(u,v)\phi^y_v(t){\rm d}v,$$
and
$$Y_u(s)=Y_1+\sum_{t\in {\cal R}}\frac{s^{\vert t\vert}}{\vert t\vert!}
\alpha(t)t!\delta_t \int_0^1 a^y(u,v)\phi^x_v(t){\rm d}v,$$
are solutions of the integral system
$$X_u(s)=X_0+s\int_0^1 a^x(u,v)\psi(Y_v(s)){\rm d}v,$$
$$Y_u(s)=Y_1+s\int_0^1 a^y(u,v)\psi(X_v(s)){\rm d}v.$$
\end{proposition}
{\it Proof:} This is consequence of Butcher's general theory (see \cite{ButcherB}).
To prove it more directly, proceed as in the proof of Theorem
\ref{DifferentialSys}

\endproof
\begin{corollary}\label{SolutionA}
Let $X_0=Y_1=0$. Assume that
$a^x(u,v)={\rm I}_{[0,u]}(v)$
and $a^y(u,v)={\rm I}_{[u,1]}(v)$.
Suppose that $\psi(z)=1/(1-z)$. Then
$$
Y_0(s)=\sum_{t\in {\cal R}}\frac{s^{\vert t\vert}}{\vert t\vert!}
\alpha(t)t!\delta_t \int_0^1 \phi^x_v(t){\rm d}v=
\sum_{t\in {\cal F}}s^{\vert t\vert}\tau(TT^*)^{\vert t\vert-1}.
$$
\end{corollary}
This result shows that the generating function of the
*-moments of the operator $T T^*$ can be obtained
by solving the system given in Proposition \ref{IntegralSystem}.
We recover in this way a result of \cite{Dykema},
Lemmas 8.5 and 8.8.
\begin{lemma}\label{DK8}
In the setting of Corollary \ref{SolutionA}, 
the generating function $Y_0(s)$ solves
\begin{equation}\label{Inverse}
G(\frac{s}{1-Y_0(s)})=s,
\end{equation}
where $G(z)=z \exp(-z)$, that is, $L(s)=s/(1-Y_0(s))$ and
$G$ are inverse with respect to composition. Moreover
$\tau(TT^*)^n=n^n/(n+1)!$.
\end{lemma}
{\it Proof:} We solve the integral system
by looking for solutions of the form
$X_u(s)=1-\exp(\lambda u)$ and
$Y_u(s)=1-\exp(\lambda'(u-1))$, with
$({\rm d}/{\rm d}u)X_u(s)=s/(1-Y_u(s))$
and $({\rm d}/{\rm d}u)Y_u(s)=-s/(1-X_u(s))$.
We deduce that $\lambda'=-\lambda$ is solution of
the equation $\lambda+s\exp(-\lambda)=0$. 
The formula for the moments of $T T^*$ is a 
consequence of Lagrange's inversion formula.

\endproof

\bigskip

\noindent {\bf Acnowledgment}

\noindent We gratefully acknowledge our debt
to the referees who carefully read the paper and
provided us with very useful comments and
relevant references.


\end{document}